\documentclass[12pt]{article}
\usepackage{fullpage,graphicx,graphics,psfrag,amsmath,amsfonts,verbatim}
\usepackage[small,bf]{caption}
\usepackage{url,ar}
\usepackage{amssymb,enumitem}

\newcommand{\BEAS}{\begin{eqnarray*}}
\newcommand{\EEAS}{\end{eqnarray*}}
\newcommand{\BEA}{\begin{eqnarray}}
\newcommand{\EEA}{\end{eqnarray}}
\newcommand{\BEQ}{\begin{equation}}
\newcommand{\EEQ}{\end{equation}}
\newcommand{\BIT}{\begin{itemize}}
\newcommand{\EIT}{\end{itemize}}
\newcommand{\BNUM}{\begin{enumerate}}
\newcommand{\ENUM}{\end{enumerate}}

\newcommand{\BA}{\begin{array}}
\newcommand{\EA}{\end{array}}

\newcommand{\ie}{{\it i.e.}}

\newcommand{\reals}{{\mbox{\bf R}}}

\newcommand{\K}{\mathcal{K}}

\makeatletter
\long\def\@makecaption#1#2{
   \vskip 9pt
   \begin{small}
   \setbox\@tempboxa\hbox{{\bf #1:} #2}
   \ifdim \wd\@tempboxa > 5.5in
        \begin{center}
        \begin{minipage}[t]{5.5in}
        \addtolength{\baselineskip}{-0.95pt}
        {\bf #1:} #2 \par
        \addtolength{\baselineskip}{0.95pt}
        \end{minipage}
        \end{center}
   \else
    \hbox to\hsize{\hfil\box\@tempboxa\hfil}
   \fi
   \end{small}\par
}
\makeatother

\newcounter{oursection}

\newcounter{lecture}

\newcounter{algorithmctr}[section]
\renewcommand{\thealgorithmctr}{\thesection.\arabic{algorithmctr}}
\newenvironment{algdesc}%
    {\refstepcounter{algorithmctr}\begin{list}{}{%
		\setlength{\rightmargin}{0\linewidth}%
		\setlength{\leftmargin}{.05\linewidth}}%
		\rmfamily\small
		\item[]{\setlength{\parskip}{0ex}\hrulefill\par%
		\nopagebreak{\bfseries\textsf{Algorithm \thealgorithmctr~}}}}%
	{{\setlength{\parskip}{-1ex}\nopagebreak\par\hrulefill} \end{list}}

\title{Disciplined Convex-Concave Programming}
\author{Xinyue Shen \and Steven Diamond \and Yuantao Gu \and Stephen Boyd}

\begin{document}
\maketitle

\begin{abstract}
In this paper we introduce \emph{disciplined convex-concave programming}
(DCCP), which combines the ideas of disciplined convex programming (DCP) with
convex-concave programming (CCP).
Convex-concave programming is an organized heuristic for solving
nonconvex problems that involve objective and constraint functions
that are a sum of a convex and a concave term.
DCP is a structured way to define convex optimization problems, based on
a family of basic convex and concave functions and a few rules
for combining them.  Problems expressed using DCP can be
automatically converted to standard form and solved by a generic solver;
widely used implementations include \texttt{YALMIP}, \texttt{CVX},
\texttt{CVXPY}, and \texttt{Convex.jl}.
In this paper we propose a framework that combines the two ideas,
and includes two improvements over previously published work on
convex-concave programming, specifically the handling of domains
of the functions, and the issue of nondifferentiability on the boundary
of the domains. We describe a Python implementation called \texttt{DCCP},
which extends \texttt{CVXPY}, and give examples.
\end{abstract}


\section{Disciplined convex-concave programming}

\subsection{Difference of convex programming}
Difference of convex (DC) programming problems have the form
\begin{equation}\label{main}
\begin{array}{ll}
\mbox{minimize} & f_0(x) - g_0(x)\\
\mbox{subject to} & f_i(x) - g_i(x) \leq 0,\quad i=1,\ldots,m,
\end{array}
\end{equation}
where $x \in \reals^n$ is the optimization variable, and the functions
$f_i: \reals^n \rightarrow \reals$ and $g_i: \reals^n \rightarrow \reals$
for $i=0,\ldots,m$ are convex.
The DC problem \eqref{main} can also include equality constraints of the
form $p_i(x)=q_i(x)$, where
$p_i$ and $q_i$ are convex; we simply
express these as the pair of inequality constraints
\[
p_i(x) - q_i(x) \leq 0, \qquad
q_i(x) - p_i(x) \leq 0,
\]
which have the difference of convex form in \eqref{main}.
When the functions $g_i$ are all affine, the problem \eqref{main} is
a convex optimization problem, and easily solved \cite{BoydandVandenberghe}.

The broad class of DC functions includes all $C^2$ functions \cite{Hartman},
so the DC problem \eqref{main} is very general.
A special case is Boolean linear programs,
which can represent many problems, such as the traveling salesman problem,
that are widely believed to be hard to solve \cite{Karp}.
DC programs arise in many applications in fields such as
signal processing \cite{Lou2015}, machine learning \cite{LeThi2008},
computer vision \cite{Lou2015Image}, and statistics \cite{Thai2014}.

DC problems can be solved globally by methods such as branch and
bound \cite{Agin, LawlerandWood}, which can be slow in practice.
Good overviews of solving DC programs globally can be found
in \cite{HorstPardalosandThoai,HorstandThoai} and the references therein.
A locally optimal (approximate) solution can be found instead through
the many techniques of general nonlinear optimization \cite{NocedalandWright}.

The convex-concave procedure (CCP) \cite{YuilleandRangarajan} is
another heuristic algorithm for finding a local optimum of \eqref{main},
which leverages our ability to efficiently solve convex optimization problems.
In its basic form,
it replaces concave terms with a convex upper bound, and then solves the resulting
convex problem, which is a restriction of the original DC problem.
Basic CCP can thus be viewed as an instance of majorization
minimization (MM)
algorithms \cite{LangeandHunter}, in which a minimization problem is
approximated by an easier to solve upper bound created around
the current point (a step called majorization) and then minimized.
Many MM extensions have been developed over the years and more can be
found in \cite{LittleandRubin, Lange2004, McLachlanandKrishnan}.
CCP can also be viewed as a version of DCA \cite{TaoandSouad}
which instead of explicitly stating the linearization, finds it by solving
a dual problem.
More information on DCA can be found at \cite{DCAweb}
and the references therein.

A recent overview of CCP, with some extensions, can be found in
\cite{Lipp2015}, where the issue of infeasibility is handled (heuristically)
by an increasing penalty on constraint violations.
The method we present in this paper is an extension of the
penalty CCP method introduced in \cite{Lipp2015}, given as
algorithm~\ref{alg:inf_start} below.

\begin{algdesc}\label{alg:inf_start}
	\emph{Penalty CCP.}
	\begin{tabbing}
		{\bf given} an initial point $x_0$, $\tau_0 > 0$, $\tau_\mathrm{max}>0$, and $\mu > 1$.\\
		$k: = 0$.\\*[\smallskipamount]
		{\bf repeat} \\
		\qquad \= 1.\ \emph{Convexify.} Form
		$\hat{g}_i(x;x_k) =g_i(x_k) + \nabla g_i(x_k)^T (x-x_k)$ for $i = 0,\ldots, m$. \\
		\> 2.\ \emph{Solve.} Set the value of $x_{k+1}$ to a solution of\\
		\> \qquad $\begin{array}{ll}
		\mbox{minimize} & f_0(x)- \hat{g}_0(x;x_k) + \tau_k\sum_{i=1}^m s_i\\
		\mbox{subject to}& f_i(x) - \hat{g}_i(x;x_k) \leq s_i,\quad i=1,\ldots,m\\
		&s_i \geq 0,\quad i = 1,\ldots,m.
		\end{array}$\\
		\> 3.\ \emph{Update $\tau$.} $\tau_{k+1} := \min(\mu\tau_k,\tau_\mathrm{max})$.\\
		\> 4.\ \emph{Update iteration.}  $k := k+1$.\\
		{\bf until} stopping criterion is satisfied.
	\end{tabbing}
\end{algdesc}
See \cite{Lipp2015} for discussion of a few variations on the penalty CCP algorithm,
such as not using slack variables for constraints that are convex, \ie,
the case when $g_i$ is affine.
Here it is assumed that $g_i$ are differentiable, and have full domain (\ie, $\reals^n$).
The first condition is not critical; we can replace $\nabla g_i(x_k)$ with a
subgradient of $g_i$ at $x_k$, if it is not differentiable.  The linearization
with a subgradient instead of the gradient is still a lower bound on $g_i$.

In some practical applications, the second assumption, that $g_i$ have full domain,
does not hold, in which case the penalty CCP algorithm can fail, by arriving
at a point $x_k$ not in the domain of $g_i$, so the convexification step fails.
This is one of the issues we will address in this paper.

\subsection{Disciplined convex programming}
Disciplined convex programming (DCP) is a methodology introduced by Grant et al.
\cite{GBY:06}
that imposes a set of conventions that must be followed when
constructing (or specifying or defining) convex programs.
Conforming problems are called \emph{disciplined convex programs}.

The conventions of DCP restrict the set of functions that can appear in
a problem and the way functions can be composed.
Every function in a disciplined convex program must come from a set
of atomic functions with known curvature and graph implementation,
or representation as partial
optimization over a cone program \cite{GrantandBoyd,NesNem:92}.
Every composition of functions $f(g_1(x),\ldots,g_k(x))$,
where $f : \reals^p \to \reals \to \reals$ is convex and
$g_1,\ldots,g_p : \reals^n \to \reals$,
must satisfy the following composition rule,
which ensures the composition is convex.
Let $\tilde{f} : \reals^p \to \reals \to \reals \cup \{\infty\}$
be the extended-value extension of $f$ \cite[Chap.~3]{BoydandVandenberghe}.
One of the following conditions must hold for each $i=1,\ldots,p$:
\begin{itemize}
\item $g_i$ is convex and $\tilde{f}$ is nondecreasing in argument $i$ on the range of $(g_1(x),\ldots,g_p(x))$.
\item $g_i$ is concave and $\tilde{f}$ is nonincreasing in argument $i$ on the range of $(g_1(x),\ldots,g_p(x))$.
\item $g_i$ is affine.
\end{itemize}
The composition rule for concave functions is analogous.
These rules allow us to certify the curvature (\ie, convexity or concavity)
of functions described as compositions using the basic atomic functions.

A DCP problem has the specific form
\begin{equation}\label{dccp_standard}
\begin{array}{ll}
\mbox{minimize/maximize} & o(x) \\
\mbox{subject to} & l_i(x) \sim r_i(x),\quad i = 1,\ldots,m,
\end{array}
\end{equation}
where $o$ (the objective), $l_i$ (lefthand sides), and $r_i$ (righthand sides)
are expressions (functions of the variable $x$) with curvature known from the
DCP rules,
and $\sim$ denotes one of the relational operators $=$, $\leq$, or $\geq$.
In DCP this problem must be convex, which imposes conditions on the curvature
of the expressions, listed below.
\BIT
\item For a minimization problem, $o$ must be convex;
for a maximization problem, $o$ must be concave.
\item When the relational operator is $=$, $l_i$ and $r_i$ must both be affine.
\item When the relational operator is $\leq$, $l_i$ must be convex, and
$r_i$ must be concave.
\item When the relational operator is $\geq$, $l_i$ must be concave, and
$r_i$ must be convex.
\EIT
Functions that are affine (\ie, are both convex and concave) can match either
curvature requirement; for example, we can minimize or maximize an affine
expression.

A disciplined convex program can be transformed into an
equivalent cone program by replacing each function with its graph implementation.
The convex optimization modeling systems
\texttt{YALMIP} \cite{Lofberg:04}, \texttt{CVX} \cite{cvx},
\texttt{CVXPY} \cite{cvxpy_paper},
and \texttt{Convex.jl} \cite{cvxjl} use DCP to verify problem convexity
and automatically convert convex programs into cone programs,
which can then be solved using generic solvers.

\subsection{Disciplined convex-concave programming}
We refer to a problem as a \emph{disciplined convex-concave program}
if it has the form \eqref{dccp_standard},
with $o$, $l_i$, and $r_i$ all having known DCP-verified curvature,
\emph{but the DCP curvature conditions for the objective and constraints need not hold}.
Such problems include DCP as a special case, but it includes many other
nonconvex problems as well.
In a DCCP problem we can, for example, maximize a convex function,
subject to nonaffine equality constraints, and nonconvex
inequality constraints between convex and concave expressions.

The general DC program~\eqref{main} and the DCCP standard form~\eqref{dccp_standard}
are equivalent.
To express \eqref{main} as \eqref{dccp_standard}, we express it as
\[
\begin{array}{ll}
\mbox{minimize} & f_0(x) - t \\
\mbox{subject to} & t = g_0(x)\\
& f_i(x) \leq g_i(x), \quad i=1,\ldots,m,
\end{array}
\]
where $x$ is the original optimization variable, and $t$ is a new optimization variable.
The objective here is convex, we have one (nonconvex) equality constraint, and
the constraints are all nonconvex (except for some special cases when $f_i$ or $g_i$
is affine)
It is straighforward to express the DCCP problem \eqref{dccp_standard} in the form
\eqref{main}, by identifying the functions $o_i$, $l_i$, and $r_i$ as $\pm f_i$ or $\pm g_i$
depending on their curvatures.

DCCP problems are an ideal standard form for DC programming because
the linearized problem in algorithm \ref{alg:inf_start} is a DCP program
whenever the original problem is DCCP.
The linearized problem can thus be automatically converted into a
cone program and solved using generic solvers.

\section{Domain and subdifferentiability}
In this section we delve deeper into an issue that is `assumed away' in the standard
treatments and discussions of DC programming, specifically, how to handle the case
when the functions $g_i$ do not have full domain.  (The functions $f_i$ can have
non-full domains, but this is handled automatically by the conversion into a cone program.)

\paragraph{An example.}
Suppose the domain of $g_i$ is $\mathcal{D}_i$, for $i=0,\ldots,m$.
If $\mathcal{D}_i \neq \reals^n$,
simply defining the linearization $\hat{g}_i(x;z)$ as
the first order Taylor expansion of $g_i$ at the point $z$ can lead to failure.
The following simple problem gives an example:
\[
\begin{array}{ll}
\mbox{minimize} & \sqrt{x} \\
\mbox{subject to} & x \geq -1,
\end{array}
\]
where $x\in \reals$ is the optimization variable.
The objective has domain $\reals_+$, and the solution is evidently $x^\star =0$.
The linearized problem in the first iteration of CCP is
\[
\begin{array}{ll}
\mbox{minimize} & x_0 + \frac{1}{2\sqrt{x_0}} (x - x_0) \\
\mbox{subject to} & x \geq -1,
\end{array}
\]
which has solution $x_1 = -1$.
The DCCP algorithm will fail in the first step of the next iteration,
since the original objective function is not defined at $x_1=-1$.

If we add the domain constraint directly into the linearized problem,
we obtain $x_1=0$,
but the first step of the next iteration also fails here, in a different way.
While $x_1$ is in the domain of the objective function, the objective is not
differentiable (or superdifferentiable) at $x_1$, so the linearization
does not exist.
This phenomenon of non-subdifferentiability or non-superdifferentiability
can only occur at a point on the boundary of the domain.

\subsection{Linearization with domain}
Suppose that the intersection of domains of all $g_i$ in problem~\eqref{main}
is $\mathcal{D} = \cap_{i=0}^m \mathcal{D}_i$.
The correct way to handle the domain is to define the linearization of $g_i$
at point $z$ to be
\begin{equation}\label{linear}
\hat{g_i}(x;z) = g_i(z) + \nabla g_i(z)^T(x - z) -  \mathcal{I}_i(x),
\end{equation}
where the indicator function is
\[
\mathcal{I}_i(x) = \left\lbrace
\begin{array}{ll}
0\;& x\in\mathcal{D}_i\\
\infty \; & x \notin \mathcal{D}_i,
\end{array}
\right.
\]
so any feasible point for the linearized problem
is in the domain $\mathcal{D}$.

Since $g_i$ is convex, $\mathcal{D}_i$ is a convex set and $\mathcal{I}_i$
is a convex function.
Therefore the `linearization' \eqref{linear} is a concave function;
it follows that if we replace the standard linearization in
algorithm \ref{alg:inf_start} with the domain-restricted linearization (\ref{linear}),
the linearized problem is still convex.

\subsection{Domain in DCCP}
Recall that we defined DCCP problems to ensure that the linearized
problem in algorithm \ref{alg:inf_start} is a DCP problem.
It is not obvious that if we replace the standard linearization with
equation (\ref{linear}) the linearized problem is still a DCP problem.
In this section we prove that the linearized DCCP problem still satisfies
the rules of DCP,
or equivalently that each $\mathcal{I}_i(x)$ has a known graph
implementation or satisfies the DCP composition rule.

If $g_i$ is an atomic function, then we assume that
\[
\mathcal{D}_i = \cap_{i=1}^p \{x \mid A_i x + b_i \in \K_i\},
\]
for some cone constraints $\K_1,\ldots,\K_p$.
The assumption is reasonable since $g_i$ itself can be represented
as partial optimization over a cone program.
The graph implementation of $\mathcal{I}_i(x)$ is simply
\[
\begin{array}{ll}
\mbox{minimize} & 0 \\
\mbox{subject to} & A_i x + b_i \in \K_i, \quad i = 1,\ldots,p.
\end{array}
\]

The other possibility is that $g_i$ is a composition of atomic functions.
Since the original problem is DCCP, we may assume that
$g_i(x) = f(h_1(x),\ldots,h_p(x))$ for some convex atomic function
$f : \reals^p \to \reals$
and DCP compliant $h_1, \ldots, h_p:\reals^n \to \reals$ such that
$f(h_1(x),\ldots,h_p(x))$ satisfies the DCP composition rule.
Then we have
\[
\mathcal{I}_i(x) = \mathcal{I}_{f}(h_1(x),\ldots,h_p(x)) +
\sum_{j=1}^p\mathcal{I}_{h_j}(x),
\]
where $\mathcal{I}_{f}$ is the indicator function for the domain of $f$
and $\mathcal{I}_{h_1},\ldots,\mathcal{I}_{h_p}$ are defined similarly.

Since $f$ is convex, $\mathcal{I}_{f}$ is convex.
Moreover, $\mathcal{I}_{f}(h_1(x),\ldots,h_p(x))$ satisfies the DCP
composition rule.
To see why, observe that for $i=1,\ldots,p$, if $h_i$ is convex
then by assumption the extended-value extension $\tilde{f}$ is
nondecreasing in argument $i$ on the range of $(h_1(x),\ldots,h_p(x))$.
It follows that $\mathcal{I}_{f}$ is nondecreasing in argument $i$
on the range of $(h_1(x),\ldots,h_p(x))$.
Similarly, if $h_i$ is concave then $\mathcal{I}_{f}$ is
nonincreasing in argument $i$ on the range of $(h_1(x),\ldots,h_p(x))$.

An inductive argument shows that
$\mathcal{I}_{h_1},\ldots,\mathcal{I}_{h_p}$ are convex and
satisfy the DCP rules.
We conclude that $\mathcal{I}_i$ satisfies the DCP composition rule.

\subsection{Sub-differentiability on boundary }
When $\mathcal{D} \neq \reals^n$,
a solution to the linearized problem $\hat{x}_k$ at iteration $k$
can be on the boundary of the closure of $\mathcal{D}$.
It is possible (as our simple example above shows)
that the convex function $g_i$ is not subdifferentiable at $\hat x_k$,
which means the linearization does not exist and the algorithm fails.
This pathology can and does occur in practical problems.

In order to handle this, at each iteration,
when the subgradient $\nabla g_i (\hat{x}_k)$ for any function $g_i$ does not exist,
we simply take a damped step,
\[
x_k = \alpha \hat{x}_k  + (1-\alpha) x_{k-1},
\]
where $0<\alpha<1$.
If $x_0$ is in the interior of the domain,
then $x_k$ will be in the interior for all $k\geq 0$,
and $\nabla g_i (x_k)$ will be guaranteed to exist.
The algorithm can (and does, for our simple example) converge to a point on the boundary
of the the domain, but each iterate is in the interior of the domain, which is enough
to guarantee that the linearization exists.

\section{Initialization}
As a heuristic method, the result of algorithm \ref{alg:inf_start} generally depends
on the initialization,
and the initial values of variables should be in the interior of the domain.
In many applications there is a natural way to carry out this initialization;
here we discuss a generic method (attempting) to do it.
Note that in general the problem of finding $x_0\in \mathcal D$ can be very hard,
so we do not expect to have a generic method that always works.

One simple and effective method is to generate
random points $z_j$ for $j=1,\ldots,k_{\mathrm{ini}}$, with entries
drawn from from i.i.d.\ standard Gaussian distributions.
We then project these points onto $\mathcal D$, \ie, solve the problems
\[
\begin{array}{ll}
\mbox{minimize} & \|x-z_j\|_2\\
\mbox{subject to} & x \in\mathcal{D},
\end{array}
\]
for $j=1, \ldots, k_\mathrm{ini}$,
denoting the solutions as $x^j_{\mathrm{ini}}$.
These points are on the boundary of $\mathcal D$ when $z_j \not\in \mathcal D$.
We then take
\[
x_0 = \frac{1}{k_{\mathrm{ini}}} \sum_{j=1}^{k_{\mathrm{ini}}}  x^j_{\mathrm{ini}}.
\]
Forming the average is a heuristic for finding $x_0$ in the interior of $\mathcal{D}$;
but it is still possible that $x_0$ is on the boundary, in which case it is an
unacceptable starting point.
As a generic practical method, however, this approach seems to work very well.

\section{Implementation}

The proposed methods described above have been implemented as the Python package
\texttt{DCCP}, which extends the package \texttt{CVXPY}.
New methods were added to \texttt{CVXPY} to return the domain of a DCP expression
(as a list of constraints), and gradients (or subgradients or supergradients) were added
to the atoms.
The linearization, damping, and initialization are handled by the package \texttt{DCCP}.
Users can form any DCCP problem of the form~\eqref{dccp_standard},
with each expression composed of functions in the \texttt{CVXPY} library.

When the \texttt{solve(method = 'dccp')} method is called on a problem object,
\texttt{DCCP} first verifies that the problem satisfies the DCCP rules.
The package then splits each non-affine equality constraint $l_i = r_i$
into $l_i \leq r_i$ and $l_i \geq r_i$.
The curvature of the objective and the left and righthand sides of each constraint
is checked, and if needed, linearized.
In the linearization the function value and gradient are \texttt{CVXPY}parameters,
which are constants whose value can change without reconstructing the problem.
For each constraint in which the left or righthand side is linearized,
a slack variable is introduced, and added to the objective.
For any expression that is linearized, the domain of the original expression
is added into the constraints.

Algorithm~\ref{alg:inf_start} is next applied to the convexified problem.
If a valid initial value of a variable is given by the user, it is used;
otherwise the generic method described above is used.
In each iteration the parameters in the linearizations (which are function and 
gradient values) are updated
based on the current value of the variables.
If a gradient (or super- or subgradient) w.r.t.\ any variable
does not exist, damping is applied to all the variables.
The convexified problem at each iteration is solved using \texttt{CVXPY}.

Some useful functions and attributes in the \texttt{DCCP} package are below.
\BIT
\item Function \verb|is_dccp(problem)| returns a boolean indicating
if an optimization problem is DCCP.
\item Attribute \verb|expression.gradient| returns a dictionary of the gradients
of a DCP expression w.r.t.\ its variables at the points specified by
 \verb|variable.value|.
 (This attribute is also in the core \texttt{CVXPY} package.)
\item Function \verb|linearize(expression)| returns the
linearization (\ref{linear}) of a DCP expression.
\item Attribute \verb|expression.domain| returns a list of constraints
describing the domain of a DCP expression.  (This attribute is also in the core
\texttt{CVXPY} package.)
\item Function \verb|convexify(constraint)| returns the transformed
constraint (without slack variables) satisfying DCP of a DCCP constraint.
\item Method \verb|problem.solve(method = 'dccp')|
carries out the proposed penalty CCP
algorithm, and returns the value of the transformed cost function,
the value of the weight $\mu_k$,
and the maximum value of slack variables at each iteration $k$.
An optional parameter is used to set the number of times to run CCP, using the
randomized initialization.
\EIT

\clearpage
\section{Examples}
In this section we describe some simple examples, show how they can be expressed
using \texttt{DCCP}, and give the results.  In each case we run the default
solve method, with no tuning or adjustment of algorithm parameters.

\subsection{Circle packing}
The aim is to arrange $n$ circles in $\reals^2$ with given radii
$r_i$ for $i=1,\ldots,n$,
so that they do not overlap and are contained in the smallest possible square \cite{Packomania}.
The optimization problem can be formulated as
\[
	\begin{array}{ll}
		\mbox{minimize} & \max_{i=1,\ldots,n} (\| c_i \|_\infty + r_i)  \\
		\mbox{subject to} & \| c_i - c_j \|_2 \geq r_i + r_j,\quad 1\leq i<j\leq n,
	\end{array}
\]
where the variables are the centers of the circles
$c_i\in\reals^2$, $i=1,\ldots,n$, and
$r_i$, $i=1,\ldots,n$, are given data.
If $l$ is the value of the objective function,
the circles are contained in the square $[-l,l]\times[-l,l]$.

This problem can be specified in \texttt{DCCP} (and solved, in the last line) as follows.
\begin{quote}
\begin{verbatim}
c = Variable(n,2)
constr = []
for i in range(n-1):
    for j in range(i+1,n):
        constr += [norm(c[i,:]-c[j,:]) >= r[i]+r[j]]
prob = Problem(Minimize(max_entries(row_norm(c,'inf')+r)), constr)
prob.solve(method = 'dccp')
\end{verbatim}
\end{quote}

\begin{figure}
	\centering
	\includegraphics[width=0.35\textwidth]{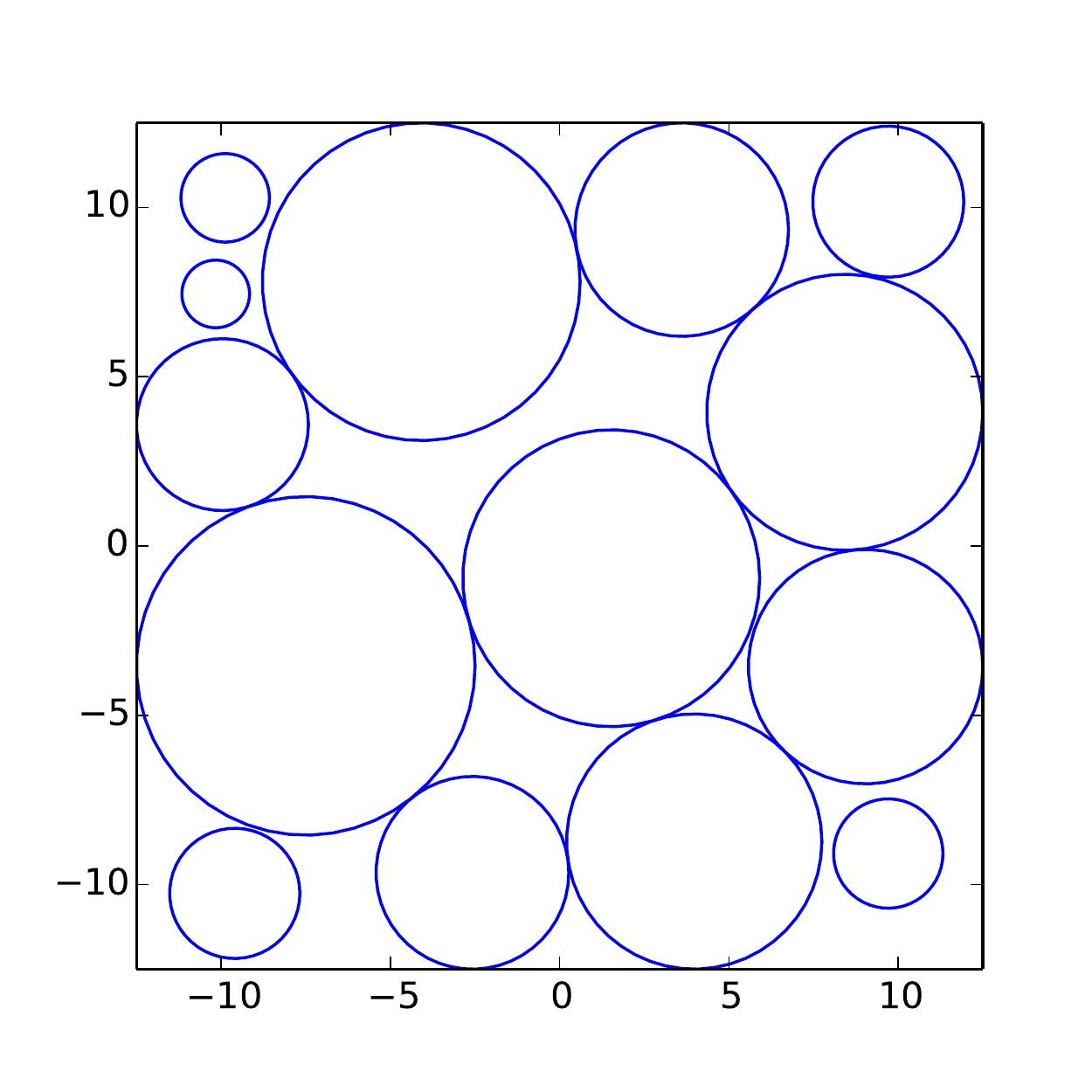}
	\caption{Circle packing.}
	\label{fig:cir_cover}
\end{figure}

The result obtained for an instance of the problem, with $n=14$ circles,
is shown in figure~\ref{fig:cir_cover}.
The fraction of the square covered by circles is $0.73$.

\clearpage
\subsection{Boolean least squares}
A binary signal $s\in\{-1,1\}^n$ is transmitted through a communication channel,
and received as $y = As + v$,
where $v\sim \mathcal{N}(0,\sigma^2 I)$ is a noise,
and $A\in\reals^{m\times n}$ is the channel matrix.
The maximum likelihood estimate of $s$ given $y$ is a solution of
\[
\begin{array}{ll}
\mbox{minimize} & \|y - Ax\|_2 \\
\mbox{subject to} & x_i^2 = 1,\quad i=1,\ldots,n,
\end{array}
\]
where $x$ is the optimization variable \cite{Forney1973}.
It is a boolean least squares problem if the objective function is squared.

The corresponding code for this problem is given below.
\begin{quote}
\begin{verbatim}
x = Variable(n)
prob = Problem(Minimize(norm(y-A*x,2)), [square(x) == 1])
result = prob.solve(method = 'dccp')
\end{verbatim}
\end{quote}
Note that the square function in the constraint is elementwise.

\begin{figure}
	\centering
	\includegraphics[width = 0.45\textwidth]{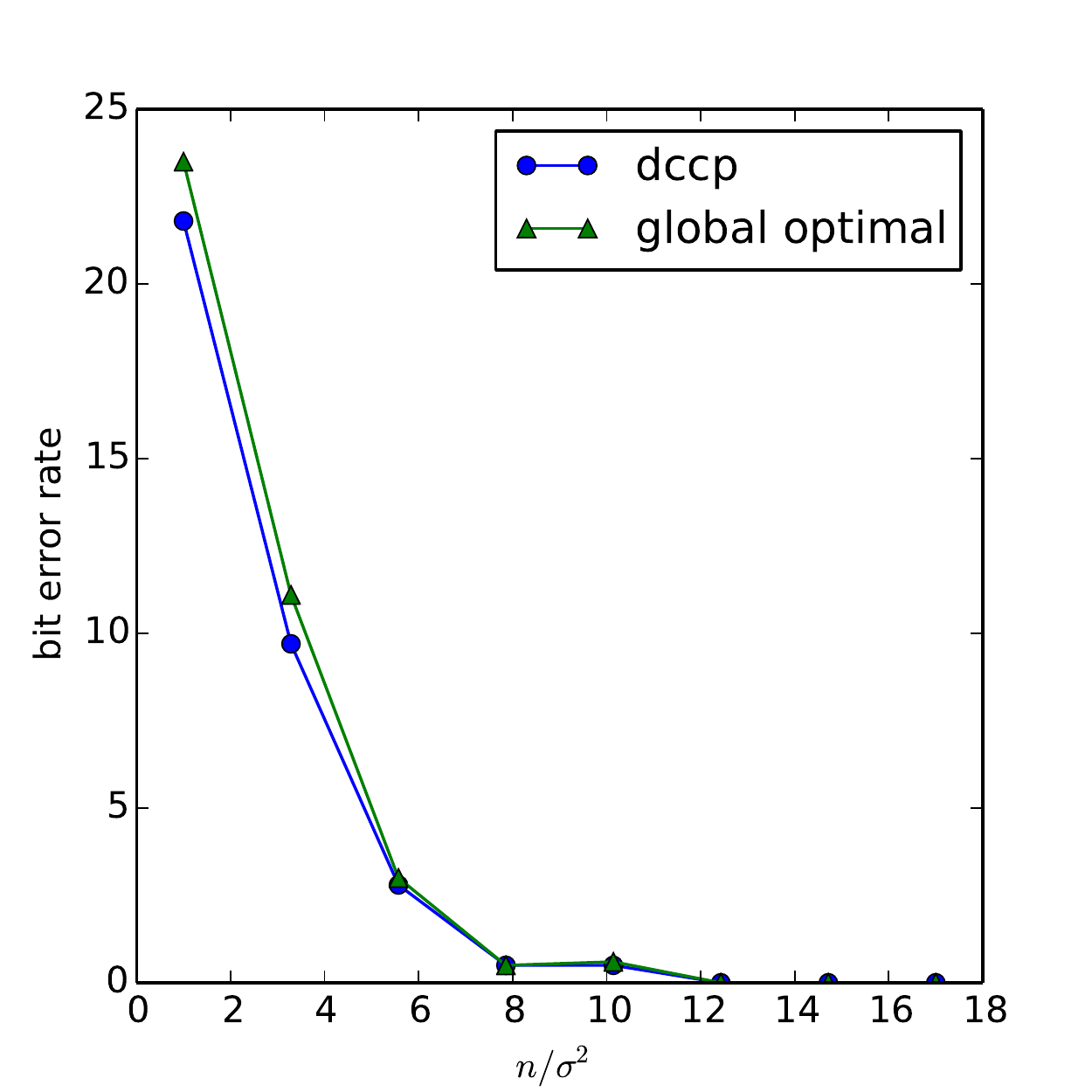}
	\caption{Boolean least squares.}
	\label{fig:binary_equa}
\end{figure}

We consider some numerical examples with
$m=n=100$, with $A_{ij}\sim \mathcal N(0,1)$ i.i.d.,
and $s_i$ i.i.d. with probability $1/2$ $1$ or $-1$.
The signal to noise ratio level is $n/\sigma^2$.
In each of the $10$ independent instances, $A$ and $s$ are generated,
and  $n/\sigma^2$ takes $8$ values from $1$ to $17$.
For each value of $n/\sigma^2$, $v$ is generated.
The bit error rates averaged from $10$ instances
are shown in figure~\ref{fig:binary_equa}.
Also shown are the same results obtained when the boolean least squares problem is
solved globally (at considerably more effort) using MOSEK \cite{mosek}.
We can see that the results, judged in terms of bit error rate, are very similar.

\clearpage
\subsection{Path planning}
The goal is to find the shortest path connecting points $a$ and $b$ in $\reals^d$
that avoids $m$ circles, centered at $p_j$
with radius $r_j$, $j = 1,\ldots,m$ \cite{Latombe:1991}.
After discretizing the arc length parametrized path into points
$x_0,\ldots,x_n$, the problem is posed as
\[
\begin{array}{ll}
\mbox{minimize} & L \\
\mbox{subject to} & x_0 = a,\quad x_n = b \\
& \|x_i - x_{i-1}\|_2 \leq L/n, \quad i=1,\ldots,n\\
& \|x_i - p_j \|_2 \geq r_j,\;\; i=1,\ldots,n,\quad j=1,\ldots,m,
\end{array}
\]
where $L$ and $x_i$ are variables,
and $a$, $b$, $p_j$, and $r_j$ are given.

The code is given below.
\begin{quote}
\begin{verbatim}
x = Variable(d,n+1)
L = Variable()
cost = L
constr = [x[:,0] == a, x[:,n] == b]
for i in range(1,n+1):
    constr += [norm(x[:,i]-x[:,i-1],2) <= L/n]
    for j in range(m):
        constr += [norm(x[:,i]-center[:,j],2) >= r[j]]
prob = Problem(Minimize(cost), constr)
result = prob.solve(method = 'dccp')
\end{verbatim}
\end{quote}

\begin{figure}
	\centering
	\includegraphics[width=0.45\textwidth]{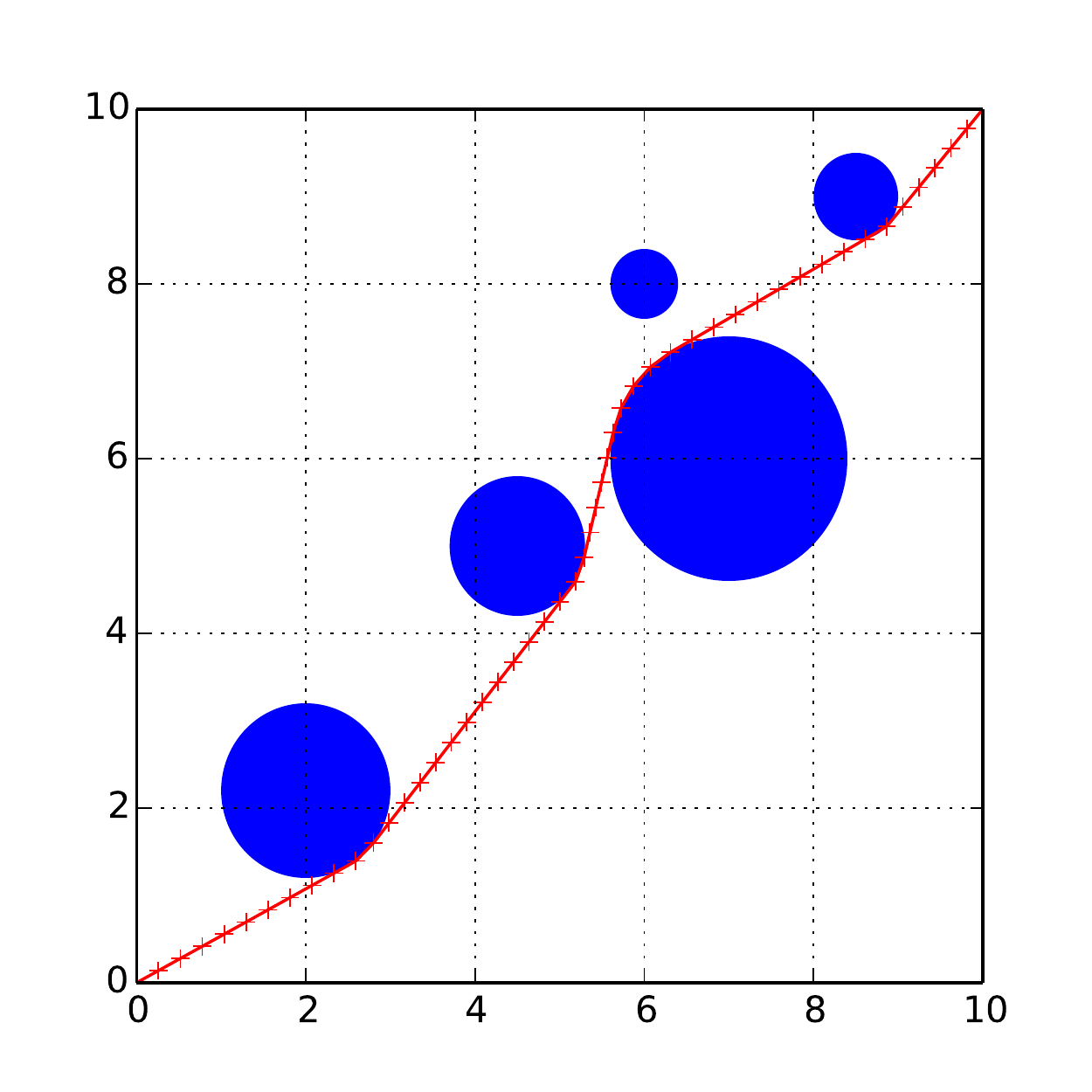}
	\caption{Path planning.}
	\label{fig:ob_avoi}
\end{figure}

An example with $d=2$ and $n=50$ is shown
in figure~\ref{fig:ob_avoi}.

\clearpage
\subsection{Control with collision avoidance}
We have $n$ linear dynamic systems, given by
\[
x^i_{t+1} = A^i x^i_t + B^i u^i_t,\quad
y^i_t = C^i x^i_t, \quad i = 1,\ldots,n,
\]
where $t=0,1,\ldots$ denotes (discrete) time, $x^i_t$ are the states, and $y^i_t$
are the outputs.
At each time $t$ for $t= 0,\ldots,T$ the $n$ outputs $y^i_t$
are required to keep a distance of at least $d_\mathrm{min}$ from each other \cite{MWCD1999}.
The initial states $x^i_0$ and ending states $x^i_n$ are given by
$x^i_\mathrm{init}$ and $x^i_\mathrm{end}$, and the inputs are
limited by $\|u^i_t\|_\infty \leq f_\mathrm{max}$.
We will minimize a sum of the $\ell_1$ norms of the inputs,
an approximation of fuel use.  (Of course we can have any convex state and input
constraints, and any convex objective.)
This gives the problem
\[
\begin{array}{ll}
\mbox{minimize} & \sum_{i=1}^n \sum_{t=0}^{T-1} \|u^i_t\|_1 \\
\mbox{subject to} & x^i_0 = x^i_\mathrm{init},\quad x^i_{T} = x^i_\mathrm{end}, \quad i=1,\ldots,n \\
& x^i_{t+1} = A^i x^i_t + B^i u^i_t, \quad t = 0,\ldots,T-1, \quad i=1,\ldots,n \\
& \|y^i_t - y^j_t\|_2 \geq d_\mathrm{min}, \quad t = 0,\ldots,T, \quad 1 \leq i <j \leq n\\
& y^i_t = C^i x^i_t,  \quad \|u^i_t\|_\infty \leq f_\mathrm{max},  \quad t = 0,\ldots,T-1, \quad i=1,\ldots,n,
\end{array}
\]
where $x^i_t$, $y^i_t$, and $u^i_t$ are variables.

The code can be written as follows.
\begin{quote}
\begin{verbatim}
constr = []
cost = 0
for i in range(n):
    for t in range(T):
        u[i] += [Variable(d)]
        constr += [norm(u[i][-1],'inf') <= f_max]
        cost += norm(u[i][-1],1)
        y[i] += [Variable(d)]
        x[i] += [Variable(2*d)]
        constr += [y[i][-1] == C[i]*x[i][-1]]
for i in range(n):
    constr += [x[i][0] == x_ini[i]]
    constr += [x[i][-1] == x_end[i]]
    for t in range(T-1):
        constr += [x[i][t+1] == A[i]*x[i][t] + B[i]*u[i][t]]
for t in range(T):
    for i in range(n-1):
        for j in range(i+1,n):
            constr += [norm(y[i][t] - y[j][t],2) >= d_min]
prob = Problem(Minimize(cost), constr)
prob.solve(method = 'dccp')
\end{verbatim}
\end{quote}

\begin{figure}
	\centering
	\includegraphics[width=\textwidth]{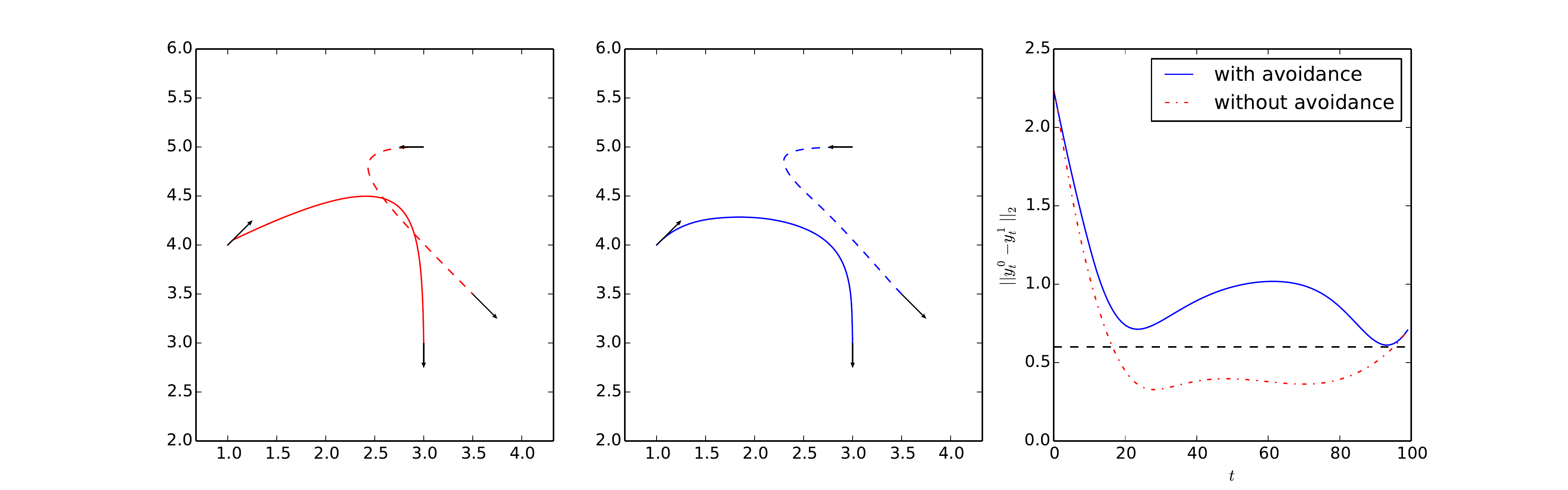}
	\caption{Optimal control with collision avoidance. \emph{Left}:
Output trajectory without collision avoidance.
\emph{Middle}: Output trajectory with collision avoidance.
\emph{Right}: Distance between outputs versus time.}
	\label{fig:colli_avoi}
\end{figure}

We consider an instance with $n=2$, with outputs (positions) $y^i_t \in \reals^2$,
$d_\mathrm{min} = 0.6$, $f_\mathrm{max} = 0.5$, $T=100$.
The linear dynamic system matrices are
\[
A^i = \begin{bmatrix}
1 & 0 & 0.1& 0\\
0 & 1 & 0 & 0.1\\
0 & 0 & 0.95 & 0\\
0 & 0 & 0 & 0.95
\end{bmatrix},\quad
B^i = \begin{bmatrix}
0 & 0\\
0 & 0\\
0.1 & 0\\
0 & 0.1
\end{bmatrix},\quad
C^i = \begin{bmatrix}
1 & 0 & 0 & 0 \\
0 & 1 & 0 & 0
\end{bmatrix}.
\]
The results are in figure~\ref{fig:colli_avoi},
where the black arrows in the first two figures show initial and final states
(position and velocity),
and the black dashed line in the third figure shows $d_\mathrm{min}$.

\clearpage
\subsection{Sparse recovery using $\ell_{1/2}$ `norm'}
The aim is to recover a sparse nonnegative signal $x_0\in\reals^n$
from a measurement vector $y = Ax_0$,
where $A\in\reals^{m\times n}$ (with $m<n$) is a known sensing matrix
\cite{Candes2008}.
A common heuristic based on convex optimization
is to minimize the $\ell_1$ norm of $x$ (which reduces here to the
sum of entries of $x$) subject to $y=Ax_0$ (and here, $x\geq 0$).
It has been proposed to minimize the sum of the \emph{squareroots} of the entries of
$x$, which since $x\geq 0$ is the same as minimizing the squareroot of the
$\ell_{1/2}$ `norm' (which is not convex, and therefore not a norm),
to obtain better recovery.
The optimization problem is
\[
\begin{array}{ll}
\mbox{minimize} & \sum_{i=1}^n \sqrt{x_i} \\
\mbox{subject to} & y = Ax,
\end{array}
\]
where $x$ is the variable. (The constraint $x\geq 0$ is implicit, since this is
the objective domain.)
This is a nonconvex problem, directly in DCCP form.

The corresponding code is as follows.
\begin{quote}
\begin{verbatim}
x = Variable(n,1)
x.value = np.ones((n,1))
prob = Problem(Minimize(sum_entries(sqrt(x))), [A*x == y])
result = prob.solve(method = 'dccp')
\end{verbatim}
\end{quote}
\begin{figure}
	\centering
	\includegraphics[width=\textwidth]{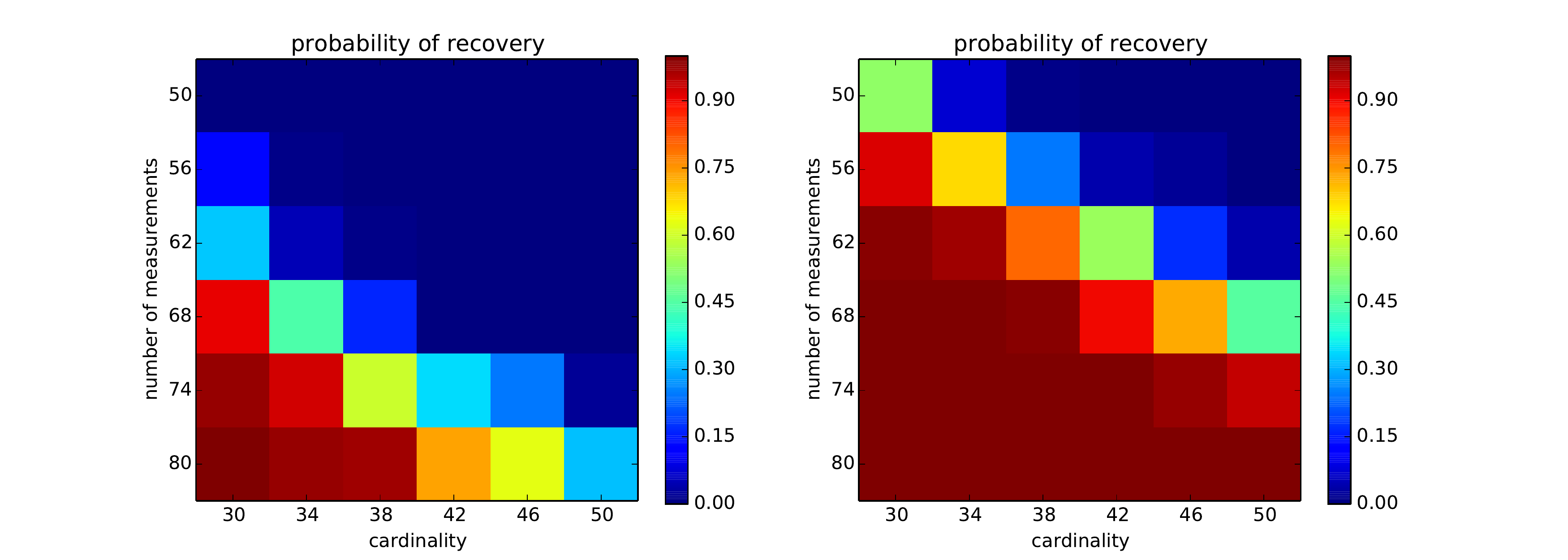}
	\caption{Sparse recovery.
	\emph{Left}: $l_1$ norm. \emph{Right}: Sqrt of $\ell_{1/2}$ `norm'.}
	\label{fig:sparse}
\end{figure}

In a numerical simulation, we take $n=100$, $A_{ij}\sim \mathcal N(0,1)$,
the positions of the nonzero entries in $x_0$ are from uniform distribution,
and the nonzero values are the absolute values of $\mathcal N(0,100)$ random variables.
To count the probability of recovery, $100$ independent instances are tested,
and a recovery is successful if the relative error $\|\hat{x}-x_0\|_2/\|x_0\|_2$
is less than $0.01$.
In each instance, the cardinality takes $6$ values from $30$ to $50$,
according to which $x_0$ is generated, and $A$ is generated for each $m$
taking one of the $6$ values from $50$ to $80$.
The results in figure~\ref{fig:sparse} verify that nonconvex recovery is
more effective than convex recovery.

\clearpage
\subsection{Phase retrieval}
Phase retrieval is the problem of recovering a signal
$x_0\in{\bf C}^n$ from the magnitudes of the complex inner products
$x_0^* a_k$, for $k=1,\ldots,m$,
where $a_k\in{\bf C}^n$ are the given measurement vectors \cite{Candes2013}.
The recovery problem can be expressed as
\[
\begin{array}{ll}
\mbox{find} & x\\
\mbox{subject to} & | x^*a_k | = y_k, \quad k=1,\ldots,m,
\end{array}
\]
where $x\in{\bf C}^n$ is the optimization variable,
and $a_k$ and $y_k\in\reals_+$ are given.
The lefthand side of the constraints are convex quadratic functions of
the real and imaginary parts of the arguments, which are in turn linear
functions of the variable $x$.

The following code segment specifies the problem.
\texttt{CVXPY} (and therefore \texttt{DCCP}) does not yet support complex
variables and constants, so
we expand complex numbers into real and imaginary parts.
\begin{quote}
\begin{verbatim}
x = Variable(2,n)
z = []
constr = []
c = np.matrix([[0,1],[-1,0]])
for k in range(m):
    z.append(Variable(2))
    z[-1].value = np.random.rand(2,1)
    constr += [norm(z[-1]) == y[k]]
    constr += [z[-1] == x*Ar[k,:] + c*x*Ai[k,:]]
prob = Problem(Minimize(0), constr)
result = prob.solve(method = 'dccp')
\end{verbatim}
\end{quote}

\begin{figure}
	\centering
	\includegraphics[width = 0.45\textwidth]{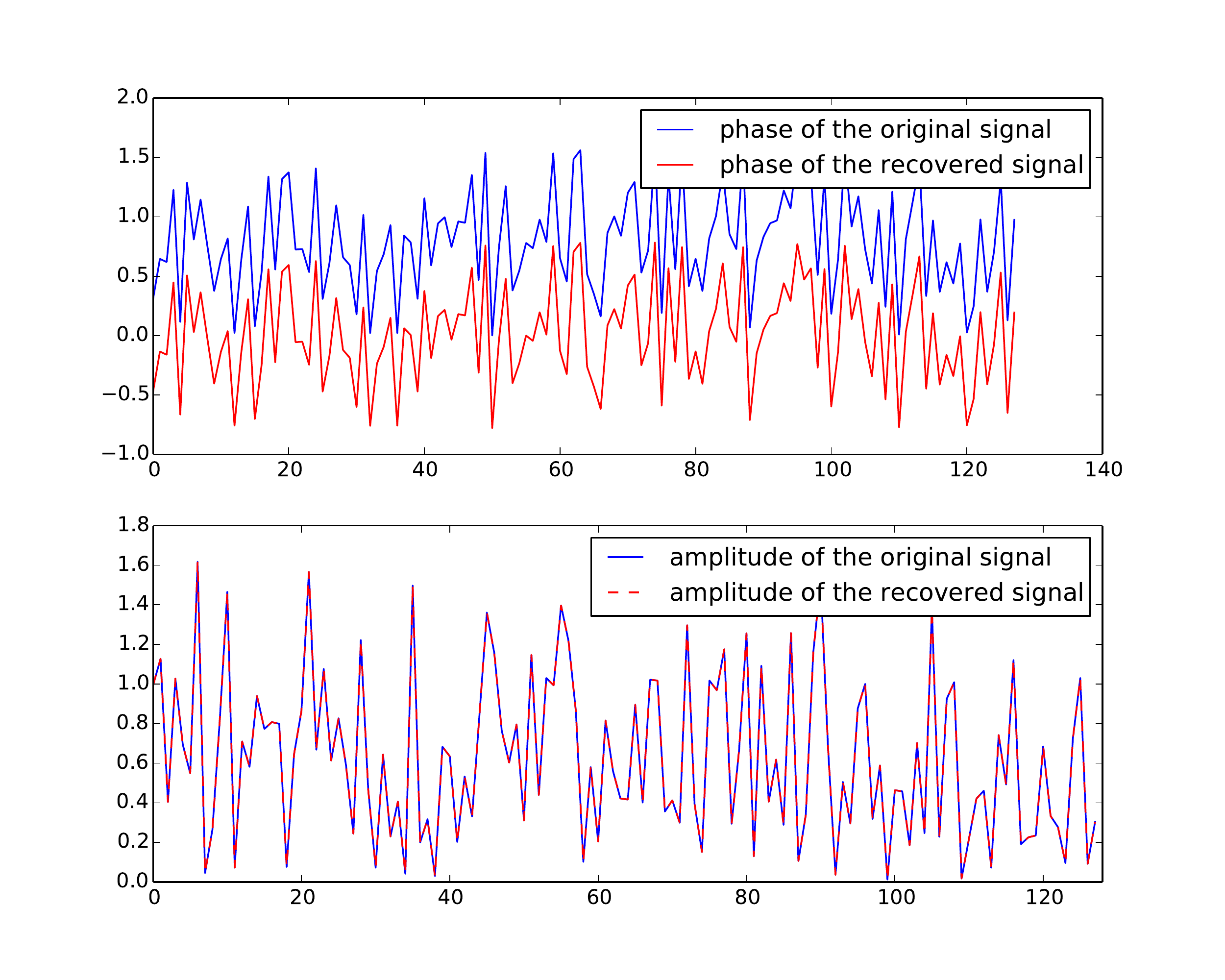}
	\caption{Phase retrieval.}
	\label{fig:phase_rec}
\end{figure}

We consider an instance with $n=128$ and $m=3n$.
The real part and the imaginary part of each entry of $x_0$ and $a_k$ are
 i.i.d.\ $\mathcal N(0,1)$.
The result in figure~\ref{fig:phase_rec} shows that
the phase is recovered (up to a global constant).

\clearpage
\subsection{Magnitude filter design}
A filter is characterized by its impulse response $\{h_k\}_{k = 1}^n$.  Its frequency
response $H:[0,\pi]\to {\bf C}$ is defined as
\[
H(\omega) = \sum_{k=1}^n h_k e^{-i\omega k},
\]
where $i=\sqrt{-1}$.
In \emph{magnitude filter design}, the goal is to find impulse response coefficients
that meet certain specifications on the magnitude of the frequency response \cite{Wu1999}.
We will consider a typical lowpass filter design problem, which can be expressed as
\[
\begin{array}{ll}
\mbox{minimize} & U_\mathrm{stop} \\
\mbox{subject to} & L_\mathrm{pass} \leq | H(\pi l/N)| \leq
U_\mathrm{pass} , \quad l=0,\ldots,l_\mathrm{pass}-1\\
&| H (\pi l/N)| \leq U_\mathrm{pass}, \quad l=l_\mathrm{pass}, \ldots,
l_\mathrm{stop}-1\\
&| H (\pi l/N)| \leq U_\mathrm{stop}, \quad l=l_\mathrm{stop}, \ldots, N,
\end{array}
\]
where $h\in\reals^n$ and $U_\mathrm{stop}\in\reals$ are the optimization variables.
The passband magnitude limits $L_\mathrm{pass}$ and $U_\mathrm{pass}$ are given.

The code can be written as follows.
\begin{quote}
\begin{verbatim}
omega = np.linspace(0,np.pi,N)
h = Variable(n)
U_stop = Variable()
constr = []
for l in range(len(omega)):
    if l < l_pass:
        constr += [norm(expo[l]*h,2) >= L_pass]
    if l < l_stop:
        constr += [norm(expo[l]*h,2) <= U_pass]
    else:
        constr += [norm(expo[l]*h,2) <= U_stop]
prob = Problem(Minimize(U_stop), constr)
result = prob.solve(method = 'dccp')
\end{verbatim}
\end{quote}

\begin{figure}
	\centering
	\includegraphics[width=0.5\textwidth]{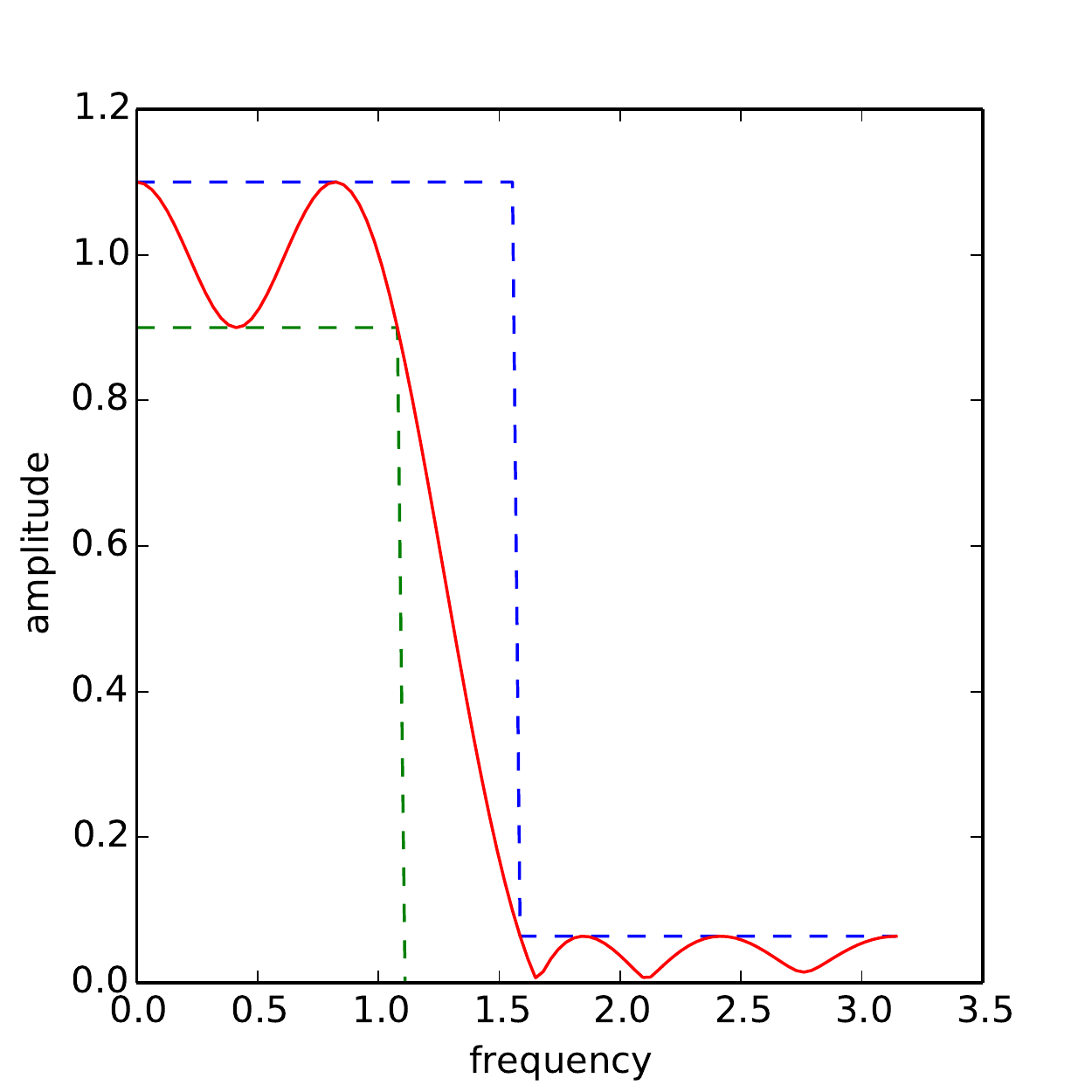}
	\caption{Low pass filter design. The frequency response magnitude upper and
lower limits are shown.}
	\label{fig:filter_design}
\end{figure}
An instance of low pass filter design, with $n=10$ and $N=100$, is shown in
figure~\ref{fig:filter_design}.

\clearpage
\subsection{Sparse singular vectors}
The left singular vectors associated with the smallest and largest singular
values of a matrix $A$ (globally) minimize and maximize $\|Ax\|_2$
subject to $\|x\|_2 = 1$.
Here we seek sparse vectors, with $\|x\|_2 =1$, which make
$\|Ax\|_2$ large or small \cite{witten2009penalized}.
To induce sparsity in $x$, we limit the $\ell_1$-norm of $x$.
(We could also limit a nonconvex sparsifier, as above in sparse recovery.)
This leads to the problems
\[
\begin{array}{ll}
\mbox{minimize/maximize} & \|Ax\|_2\\
\mbox{subject to} & \|x\|_2 = 1, \quad \|x\|_1 \leq \mu,
\end{array}
\]
where $x\in\reals^n$ is the variable and $\mu\geq 0$ controls the sparsification,
to find $x$ that is sparse, satisfies $\|x\|_2=1$,
and makes $\|Ax\|_2$ small or large.  We call such a vector, with some abuse of notation,
a \emph{sparse singular vector}.
Since $\|x\|_2=1$, we know $1\leq \|x\|_1 \leq \sqrt{n}$, so the range of $\mu$
can be set as $[1,\sqrt{n}]$.

The code (for minimization) is the following.
\begin{quote}
	\begin{verbatim}
	x = Variable(n)
	prob = Problem(Minimize(norm(A*x)), [norm(x) == 1, norm(x,1) <= mu])
	prob.solve(method = 'dccp')
	\end{verbatim}
\end{quote}

\begin{figure}
	\centering
	\includegraphics[width=0.6\textwidth]{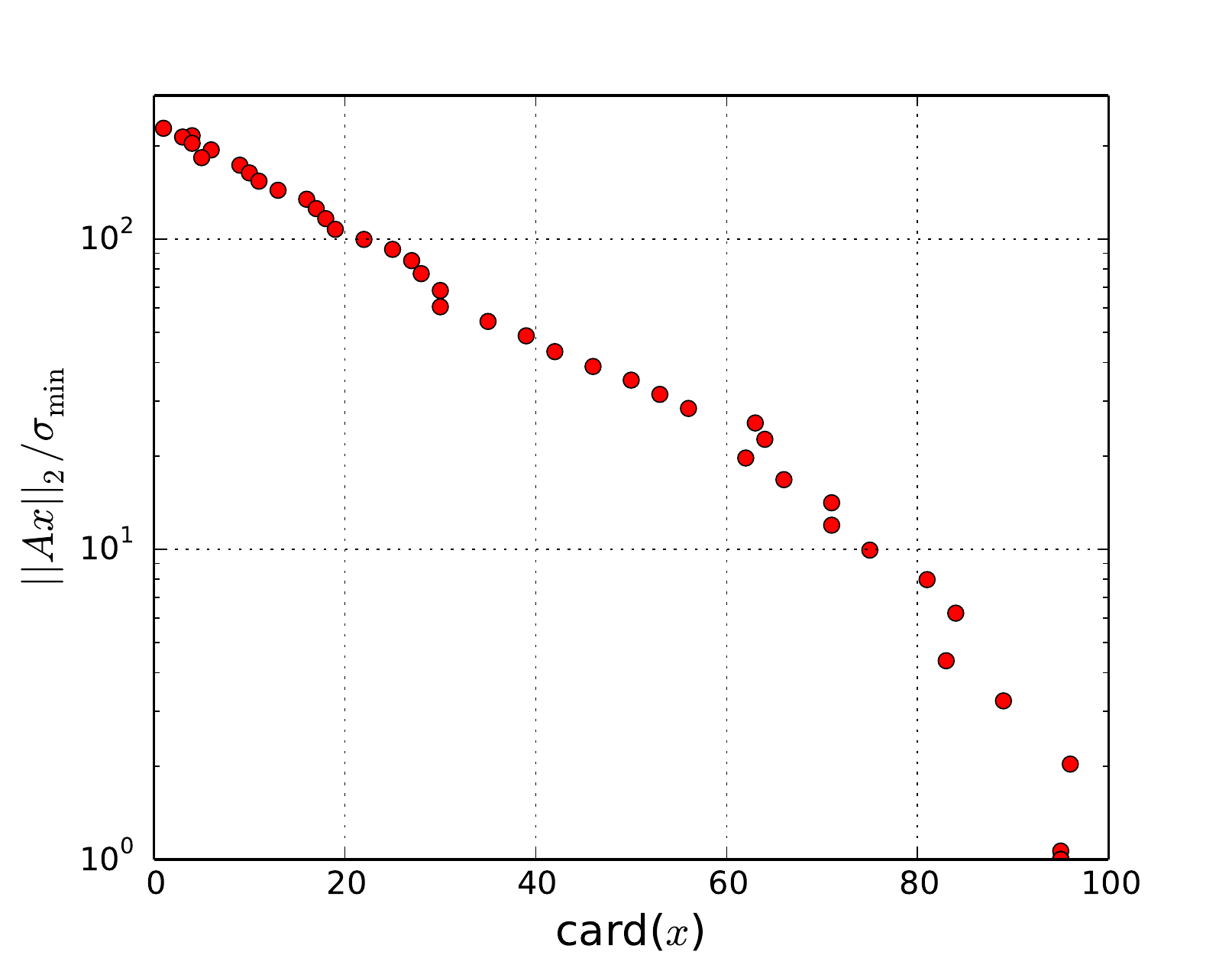}
	\caption{Sparse singular vectors.}
	\label{fig:sparse_minsingular}
\end{figure}

We consider an instance for minimization with a random matrix 
$A\in\reals^{100\times 100}$ with i.i.d.\ entries $A_{ij}\sim\mathcal{N}(0,1)$,
with (positive) smallest singular value $\sigma_{\mathrm{min}}$.
The parameter $\mu$ is swept from $1$ to $10$ with increment $0.2$,
and for each value of $\mu$ the result of solving the problem above 
is shown as a red dot in figure~\ref{fig:sparse_minsingular}.
The most left point in the figure corresponds to $\|x\|_1\leq1$,
which gives cardinality $1$.  (In this instance it achieves the globally optimal
value, which is the smallest of the norm of the columns of $A$.) 

\clearpage
\subsection{Gaussian covariance matrix estimation}
Suppose $y_i\in\reals^n$ for $i = 1,\ldots,N$ are points drawn i.i.d\
from $\mathcal{N}(0,\Sigma)$.
Our goal is to estimate the parameter $\Sigma$ given these samples.
The maximum likelihood problem of estimating $\Sigma$ is convex in
the inverse of $\Sigma$, but not $\Sigma$ \cite{Banerjee:2006}.
If there are no other constraints on $\Sigma$, the maximum likelihood estimate
is $\hat \Sigma = \frac{1}{N}\sum_{i=1}^N y_i y_i^T$, the empirical covariance
matrix.
We consider here the case where
the sign of the off-diagonal entries in $\Sigma$ is known; that is, we know which
entries of $\Sigma$ are negative, which are zero, and which are positive.
(So we know which components of $y$ are uncorrelated, and which are negatively
and positively correlated.)

The maximum likelihood problem is then
\[
\begin{array}{ll}
\mbox{maximize} & -\log\det (\Sigma) -\frac{1}{N}\sum_{i=1}^N y_i^T \Sigma^{-1} y_i \\
\mbox{subject to} & \Sigma_{\Omega_+} \geq 0, \quad \Sigma_{\Omega_-} \leq 0, \quad \Sigma_{\Omega_0} = 0,
\end{array}
\]
where $\Sigma$ is the variable,
and the index sets $\Omega_+$, $\Omega_-$, and $\Omega_0$ are given.
The objective is a difference of convex functions, so we
transform the problem into the following DCCP problem with additional variable $t$,
\[
\begin{array}{ll}
\mbox{maximize} & -\log\det (\Sigma) - t \\
\mbox{subject to} & \frac{1}{N}\sum_{i=1}^N y_i^T \Sigma^{-1} y_i \leq t\\
& \Sigma_{\Omega_+} \geq 0, \quad \Sigma_{\Omega_-} \leq 0, \quad \Sigma_{\Omega_0} = 0.
\end{array}
\]

The code is as follows.
\begin{quote}
\begin{verbatim}
Sigma = Variable(n,n)
t = Variable()
cost = -log_det(Sigma) - t
trace_val = trace(sum([matrix_frac(y[:,i], Sigma)/N for i in range(N)]))
prob = Problem(Maximize(cost),
               [trace_val <= t,
                Sigma[pos] >= 0,
               	Sigma[neg] <= 0,
               	Sigma[zero] == 0])
prob.solve(method = 'dccp')
\end{verbatim}
\end{quote}

\begin{figure}
	\centering
	\includegraphics[width=\textwidth]{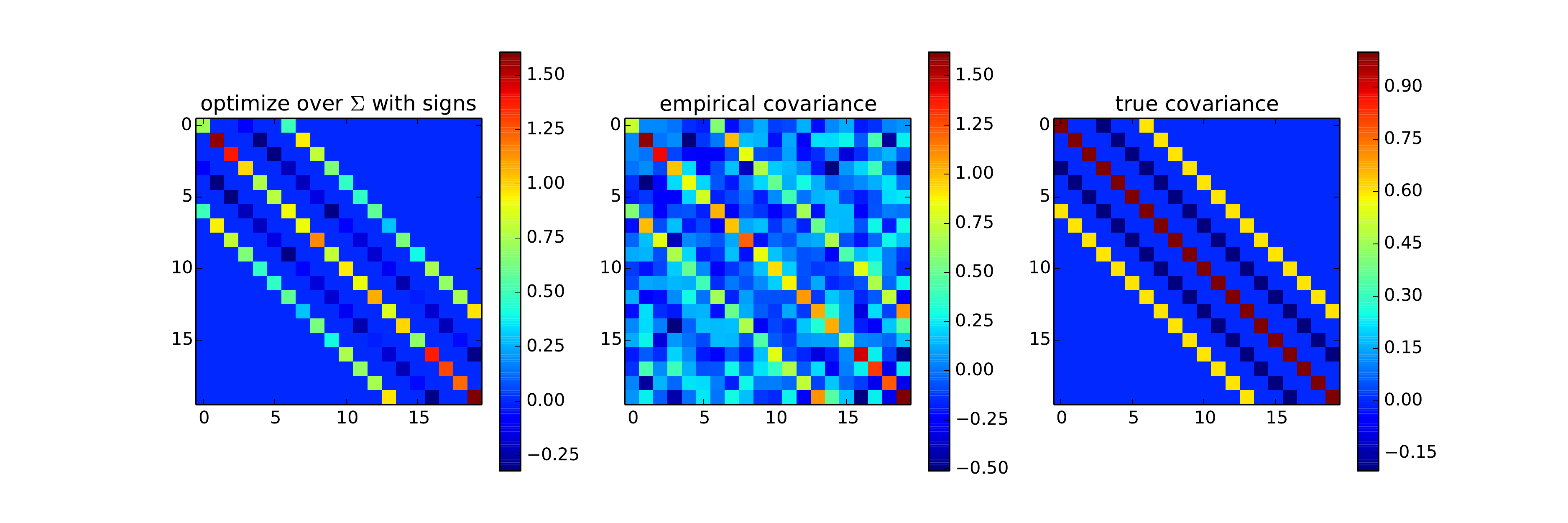}
	\caption{Gaussian covariance matrix estimation.}
	\label{fig:cov}
\end{figure}

An example with $n=20$ and $N=30$ is in figure~\ref{fig:cov}.   Not surprisingly,
knowledge of the signs of the entries of $\Sigma$ allows us to obtain a much
better estimate of the true covariance matrix.

\clearpage
\section*{Acknowledgments}
This material is based
upon work supported by the National Science Foundation Graduate Research Fellowship
under Grant No. DGE-114747, by the DARPA X-DATA and SIMPLEX programs,
and by the CSC State Scholarship  Fund.

\bibliography{dccp}

\newcommand{\etalchar}[1]{$^{#1}$}
\begin{thebibliography}{MWCD99}

\bibitem[Agi66]{Agin}
N.~Agin.
\newblock Optimum seeking with branch and bound.
\newblock {\em Management Science}, 13:176--185, 1966.

\bibitem[ALNT08]{LeThi2008}
L.~T.~H. An, H.~M. Le, V.~V. Nguyen, and P.~D. Tao.
\newblock A {DC} programming approach for feature selection in support vector
  machines learning.
\newblock {\em Advances in Data Analysis and Classification}, 2(3):259--278,
  2008.

\bibitem[An15]{DCAweb}
L.~T.~H. An.
\newblock {DC} programming and {DCA}: local and global approaches - theory,
  algorithms and applications.
\newblock \url{http://lita.sciences.univ-metz.fr/~lethi/DCA.html}, 2015.

\bibitem[ApS15]{mosek}
MOSEK ApS.
\newblock {\em The MOSEK Python optimizer API manual. Version 7.1 (Revision
  49).}, 2015.

\bibitem[BGdN06]{Banerjee:2006}
O.~Banerjee, L.~E. Ghaoui, A.~d'Aspremont, and G.~Natsoulis.
\newblock Convex optimization techniques for fitting sparse gaussian graphical
  models.
\newblock In {\em Proceedings of the 23rd International Conference on Machine
  Learning}, ICML '06, pages 89--96, 2006.

\bibitem[BV04]{BoydandVandenberghe}
S.~Boyd and L.~Vandenberghe.
\newblock {\em Convex Optimization}.
\newblock Cambridge University Press, Cambridge, 2004.

\bibitem[CESV13]{Candes2013}
E.~J. Cand{\`{e}}s, Y.~C. Eldar, T.~Strohmer, and V.~Voroninski.
\newblock Phase retrieval via matrix completion.
\newblock {\em SIAM Journal on Imaging Sciences}, 6(1):199--225, 2013.

\bibitem[{CVX}12]{cvx}
{CVX Research, Inc.}
\newblock {CVX}: Matlab software for disciplined convex programming, version
  2.0.
\newblock \url{http://cvxr.com.cvx}, August 2012.

\bibitem[CW08]{Candes2008}
E.~J. Cand{\`{e}}s and M.~B. Wakin.
\newblock An introduction to compressive sampling.
\newblock {\em IEEE Signal Processing Magazine}, 25(2):21--30, 2008.

\bibitem[DB16]{cvxpy_paper}
S.~Diamond and S.~Boyd.
\newblock {CVXPY}: A {P}ython-embedded modeling language for convex
  optimization.
\newblock {\em To appear, Journal of Machine Learning Research}, 2016.

\bibitem[For72]{Forney1973}
G.~Forney.
\newblock Maximum-likelihood sequence estimation of digital sequences in the
  presence of intersymbol interference.
\newblock {\em IEEE Transactions on Information Theory}, 18(3):363--378, 1972.

\bibitem[GB08]{GrantandBoyd}
M.~Grant and S.~Boyd.
\newblock Graph implementation for nonsmooth convex programs.
\newblock In V.~Blondel, S.~Boyd, and H.~Kimura, editors, {\em Recent Advances
  in Learning and Control}, Lecture Notes in Control and Information Sciences.
  Springer-Verlag Limited, 2008.
\newblock \url{http://stanford.edu/~boyd/graph_dcp.html}.

\bibitem[GBY06]{GBY:06}
M.~Grant, S.~Boyd, and Y.~Ye.
\newblock Disciplined convex programming.
\newblock In L.~Liberti and N.~Maculan, editors, {\em Global Optimization: From
  Theory to Implementation}, Nonconvex Optimization and its Applications, pages
  155--210. Springer, 2006.

\bibitem[Har59]{Hartman}
P.~Hartman.
\newblock On functions representable as a difference of convex functions.
\newblock {\em Pacific Journal of Math}, 9(3):707--713, 1959.

\bibitem[HPT95]{HorstPardalosandThoai}
R.~Horst, P.~M. Pardalos, and N.~V. Thoai.
\newblock {\em Introduction to Global Optimization}.
\newblock Kluwer Academic Publishers, Dordrecht, Netherlands, 1995.

\bibitem[HT99]{HorstandThoai}
R.~Horst and N.~V. Thoai.
\newblock {DC} programming: overview.
\newblock {\em Journal of Optimization Theory and Applications}, 103(1):1--43,
  1999.

\bibitem[Kar72]{Karp}
R.~M. Karp.
\newblock Reducibility among combinatorial problems.
\newblock In R.~E. Miller and J.~W. Thatcher, editors, {\em Complexity of
  Computer Computation}, pages 85--104. Plenum, 1972.

\bibitem[Lan04]{Lange2004}
K.~Lange.
\newblock {\em Optimization}.
\newblock Springer Texts in Statistics. Springer, New York, New York, 2004.

\bibitem[Lat91]{Latombe:1991}
J.~C. Latombe.
\newblock {\em Robot Motion Planning}.
\newblock Kluwer Academic Publishers, 1991.

\bibitem[LB15]{Lipp2015}
T.~Lipp and S.~Boyd.
\newblock Variations and extension of the convex--concave procedure.
\newblock {\em Optimization and Engineering}, pages 1--25, 2015.

\bibitem[LHY00]{LangeandHunter}
K.~Lange, D.~R. Hunter, and I.~Yang.
\newblock Optimization transfer using surrogate objective functions.
\newblock {\em Journal of Computational and Graphical Statistics}, 9(1):1--20,
  2000.

\bibitem[Lof04]{Lofberg:04}
J.~Lofberg.
\newblock {YALMIP}: A toolbox for modeling and optimization in {MATLAB}.
\newblock In {\em Proceedings of the {IEEE} International Symposium on Computed
  Aided Control Systems Design}, pages 294--289, September 2004.

\bibitem[LOX15]{Lou2015}
Y.~Lou, S.~Osher, and J.~Xin.
\newblock {\em Computational aspects of constrained L1-L2 minimization for
  compressive sensing}, volume 359 of {\em Advances in Intelligent Systems and
  Computing}, pages 169--180.
\newblock 2015.

\bibitem[LR87]{LittleandRubin}
R.~J.~A. Little and D.~B. Rubin.
\newblock {\em Statistical Analysis with Missing Data}.
\newblock John Wiley \& Sons, New York, New York, 1987.

\bibitem[LW66]{LawlerandWood}
E.~L. Lawler and D.~E. Wood.
\newblock Branch-and-bound methods: a survey.
\newblock {\em Operations Research}, 14:699--719, 1966.

\bibitem[LZOX15]{Lou2015Image}
Y.~Lou, T.~Zeng, S.~Osher, and J.~Xin.
\newblock A weighted difference of anisotropic and isotropic total variation
  model for image processing.
\newblock {\em SIAM Journal on Imaging Sciences}, 8(3):1798--1823, 2015.

\bibitem[MK07]{McLachlanandKrishnan}
G.~McLachlan and T.~Krishnan.
\newblock {\em The EM algorithm and extensions}.
\newblock John Wiley \& Sons, 2007.

\bibitem[MWCD99]{MWCD1999}
A.~Miele, T.~Wang, C.~S. Chao, and J.~B. Dabney.
\newblock Optimal control of a ship for collision avoidance maneuvers.
\newblock {\em Journal of Optimization Theory and Applications},
  103(3):495--519, 1999.

\bibitem[NN92]{NesNem:92}
Y.~Nesterov and A.~Nemirovsky.
\newblock Conic formulation of a convex programming problem and duality.
\newblock {\em Optimization Methods and Software}, 1(2):95--115, 1992.

\bibitem[NW06]{NocedalandWright}
J.~Nocedal and S.~J. Wright.
\newblock {\em Numerical Optimization}.
\newblock Springer, 2006.

\bibitem[Spe13]{Packomania}
E.~Specht.
\newblock Packomania.
\newblock \url{http://www.packomania.com/}, October 2013.

\bibitem[THA{\etalchar{+}}14]{Thai2014}
J.~Thai, T.~Hunter, A.~K. Akametalu, C.~J. Tomlin, and A.~M. Bayen.
\newblock Inverse covariance estimation from data with missing values using the
  concave-convex procedure.
\newblock In {\em Decision and Control (CDC), 2014 IEEE 53rd Annual Conference
  on}, pages 5736--5742. IEEE, 2014.

\bibitem[TS86]{TaoandSouad}
P.~D. Tao and E.~B. Souad.
\newblock Algorithms for solving a class of nonconvex optimization problems.
  {M}ethods of subgradients.
\newblock In J.-B. Hiriart-Urruty, editor, {\em {FERMAT} Days 85: Mathematics
  for Optimization}, pages 249--271. Elsevier Scince Publishers B. V., 1986.

\bibitem[UMZ{\etalchar{+}}14]{cvxjl}
M.~Udell, K.~Mohan, D.~Zeng, J.~Hong, S.~Diamond, and S.~Boyd.
\newblock Convex optimization in {J}ulia.
\newblock In {\em Proceedings of the Workshop for High Performance Technical
  Computing in Dynamic Languages}, pages 18--28, 2014.

\bibitem[WBV99]{Wu1999}
S.~P. Wu, S.~Boyd, and L.~Vandenberghe.
\newblock {\em Applied and Computational Control, Signals, and Circuits: Volume
  1}, chapter FIR Filter Design via Spectral Factorization and Convex
  Optimization, pages 215--245.
\newblock Birkh{\"a}user Boston, 1999.

\bibitem[WTH09]{witten2009penalized}
D.~M. Witten, R.~Tibshirani, and T.~Hastie.
\newblock A penalized matrix decomposition, with applications to sparse
  principal components and canonical correlation analysis.
\newblock {\em Biostatistics}, 2009.

\bibitem[YR03]{YuilleandRangarajan}
A.~L. Yuille and A.~Rangarajan.
\newblock The concave-convex procedure.
\newblock {\em Neural Computation}, 15(4):915--936, 2003.

\end{thebibliography}

\end{document}